\documentclass[11pt]{amsart}
\usepackage{txfonts}
\usepackage{}
\usepackage{amssymb}
\theoremstyle{plain}
\newtheorem{Thm}{Theorem}

\newtheorem{Def}[Thm]{Definition}

\begin{document}

\title[Discrete Morse flow for Ricci flow]
{Discrete Morse flow for Ricci flow and Porous Media equation}

\author{Li Ma, Ingo Witt}
\address{Distinguished Professor, Department of mathematics \\
Henan Normal university \\
Xinxiang, 453007 \\
China} \email{lma@tsinghua.edu.cn}

\address{Ingo Witt, Math. Institut \\
Universitat G\"ottingen \\
Bunsenstr. 3-5, D-37073, G\"ottingen, Germany}

\email{iwitt@uni-math.gwdg.de}

\thanks{The research is partially supported by the National Natural Science
Foundation of China 10631020 and SRFDP 20090002110019}

\begin{abstract}
In this paper, we study the discrete Morse flow for the Ricci flow
on football, which is the 2-sphere with removed north and south
poles and with the metric $g_0$ of constant scalar curvature, and
and for Porous media equation on a bounded regular domain in the
plane. We show that with a suitable assumption about $g(0)$ we have
a weak approximated discrete Morse flow for the approximated Ricci
flow and Porous media equation on any time intervals.


\end{abstract}

 \maketitle

\section{Introduction}
There are relative few results about computational models for the
Ricci flow in two dimensions. The purpose of this paper is to try
this area by giving some approximated computational models, namely
the discrete Morse flow for the 2-d Ricci flow. We shall first
consider the Porous-media equation on a bounded regular domain in
the plane. As is well-known that the limiting equation of the
Porous-Media equation is the Ricci flow on the domain. Then we may
consider this discrete flow as the approximated computational scheme
for the Ricci flow. We also consider two different modelings of the
Ricci flow on the singular surface, the American football. Our
method can also be worked out in a similar way in the regular
spherical surfaces and other surfaces.

Since the geometric and analytic parts of singular surfaces are not
well-known, we now recall the geometry of a special singular
surface, namely, the American football. Let $S=S_0$ be the sphere
with north and south poles removed and with the metric
$$
g_0=dr^2+(\alpha \sin r)^2d\theta^2,
$$
where $\alpha\in (0,1)$, $0\leq r\leq \pi$ such that $r=0$
corresponding the north pole,  and $0\leq \theta\leq 2\pi$. With
this metric, we can see that $S_0$ has two singularities of equal
angle $2\pi\alpha$ at the north and south poles. By a direct
computation we know that the scalar curvature of the metric $g_0$ is
equal to $2$. Recall here that the Laplacian operator of $g_0$ on
function is
$$
\Delta_0=\partial_r^2+\frac{\cos r}{\sin
r}\partial_r+\frac{1}{(\alpha \sin r)^2}\partial^2_{\theta}
$$
and the area element is
$$
dA=\alpha \sin r drd\theta.
$$
The area of $(S,g_0)$ is $|A|=4\pi\alpha$. The aim of this note is
to study the normalized Ricci flow on $(S,g_0)$:
$$
\partial_t g=(\rho-R)g
$$
where $g=g(t)=e^{2u(t)}g_0$, $R=R(g)$ is the scalar curvature of the
metric $g$, and $\rho$ is some real constant (and we may choose it
such that the area of the flow $g(t)$ is constant). We may write the
evolution equation of the Ricci flow as
$$
e^u\partial_t e^u=\frac{1}{2}\rho e^{2u}+\Delta u-1,
$$
where $\Delta$ is the Laplacian operator of the metric $g_0$. The
standard example is that for the football metric
$g_0=dr^2+(\alpha\sin r)^2d\theta^2$ $(0<\alpha<1$) the
(un-normalized) Ricci flow is $g(t)=(1-2t)g_0$ for $t<\frac{1}{2}$,
which extinct at $t=\frac 1 2$.

For the construction of the discrete Morse flow for the 2-d Ricci
flow, we shall rely on some variational structure related. Hence we
need some analytical part of the related variational functionals.
 Here we recall one important
consequence of a Moser type inequality proven by W.Chen and C.Li
\cite{CL}. For any $\phi\in H^1(S)$,
\begin{equation}\label{chen}
\int_S e^\phi dA\leq C_2exp\{\frac{1}{16\pi\alpha}\int_S|\nabla
\phi|^2dA+\frac{1}{|S|}\int \phi dA\}
\end{equation}
and the functional
$$
\phi\in H^1(S)\to \int_S e^\phi dA
$$
is continuous with respect to the weak convergence in $H^1$. For our
purpose, we reformulate the inequality above as For any $u\in
H^1(S)$,
\begin{equation}\label{chen1}
\fint_S e^{2u} dA\leq C_2exp\{\fint_S|\nabla u|^2dA+\fint_S 2u dA\}
\end{equation}
with $\fint udA=\frac{1}{|S|}\int u dA$. We shall use this fact to
get a weak solution to the Ricci flow on $(S,g_0)$.

The precise results will be stated and proved in the following three
sections. Similar results for Yamabe flow and higher dimensional
Porous-media equations are also true. The plan of this paper is
below. In section \ref{porous}, we study the existence of a weak
solution to the porous media equation in a bounded regular domain.
In section \ref{sect3}, we introduce the method of the discrete
Morse flow method to an approximated ricci flow in the singular
surface with symmetry the initial data. We study a perturbated Ricci
flow model in the last section.

\section{discrete Morse flow for the Porous-media
equation}\label{porous}

 The existence of the weak solution to the
Porous-media equation can be done by using the Galerkin method
\cite{Ar}. Other related methods can be found in \cite{DK}. Here we
propose a new method, which is the adapted discrete Morse flow
method \cite{K}.

  To make our idea more clear, we start from the discrete Morse flow
  for the Porous-Media equation on the bounded regular domain $\Omega$
  in the plane $R^2$. In some sense our domain be be any one with
  the Sobolev imbedding theorem holds true.

  Assume that $m>0$. Given an initial regular data $u_0$ and any $T>0$. We consider the porous media
  equation
  \begin{equation}\label{pm}
\partial_t u=\Delta u^m, \ \ \  \Omega\times [0,T], \ \ \
  \end{equation}
  with the initial data $u|_{t=0}=u_0$ and with the boundary condition $u(t)=u_0$ on $\partial \Omega\times t$.

This equation has a very close relation with the Ricci flow. In
fact, by taking the limit $m\to 0$, the limiting equation of
(\ref{pm}) is
\begin{equation}\label{pm1}
\partial_t u=\Delta \log u,
\end{equation}
which is the Ricci flow on the plane $R^2$.

To introduce the discrete Morse flow for (\ref{pm}) we set $v=u^m$
and let $\alpha=1/m$. Then (\ref{pm}) is induced into
\begin{equation}\label{pm2}
\partial_t v^\alpha=\Delta v
\end{equation}
with the initial data $v=v_0=u_0^\alpha$ at $t=0$. We assume that
$v_0\in H^1(\Omega)$. The key idea is that we set $\alpha=2\beta-1$
with $\beta>1/2$ and make $$\partial_t v^\alpha= C_\beta
v^{\beta-1}\partial_t v^\beta,$$ where $C_\beta>0$ such that
$\frac{\beta}{2\beta-1}C_\beta=1$.

In below, we assume that $\beta>1$ and make the notation that
$v^\beta=|v|^\beta$ for any $v$. For any $N>1$ be a large integer
and for any $T>0$,let
$$
h=T/N, \ \ t_n=nh, \ \ n=0,1,2,...,N.
$$

Assume that we have constructed $v_j\in H^1(S)$, $0\leq j\leq n-1$
and $u_{n-1}$ is a minimizer of the functional
$$
I_{n-1}(u)=\frac{C_\beta}{2h}\int_\Omega
|v^\beta-v_{n-2}^\beta|^2dx+\frac{1}{2}\int_\Omega |\nabla v|^2dx
$$
 on $H:=\{v\in H^1(\Omega); v-v_0\in H_0^1(\Omega)\}$. Define
$$
I_n(u)=\frac{C_\beta}{2h}\int_\Omega
|v^\beta-v_{n-1}^\beta|^2dx+\frac{1}{2}\int_\Omega |\nabla v|^2dx
$$
on $H$. It is clear that the infimum is finite and by using the
Poincare inequality to $v-v_0$, any minimizing sequence is bounded
in $H$. We remark here that we may choose the functions in the
minimizing sequence non-negative since our functional is even
functional \cite{A}. By the direct method, we know that $I_n$ has an
unique minimizer $u_n$ in $H$, which satisfies that
$$
\frac{\beta C_\beta}{h}(v^\beta-v_{n-1}^\beta)v^{\beta-1}=\Delta v
$$
with the uniform energy bound
$$\frac{C_\beta}{2h}\int_\Omega
|v_n^\beta-v_{n-1}^\beta|^2dx+\frac{1}{2}\int_\Omega |\nabla
v_n|^2dx\leq \frac{1}{2}\int_\Omega |\nabla v_{n-1}|^2dx\leq C.
$$

We inductively define $v_N(t)\in H^1$ for $t\in [-h,T]$ such that
for $n=1,...,N$,
$$
v_N(t)=v_n,
$$
on $[t_{n-1},t_n]$ and $v_N(t)=v_0$ on $[-h,0]$. We define, for
$1\leq n\leq N$,
$$
\partial_t v_N^{\beta}(t)=\frac{1}{h}(v_n^\beta-u_{n-1}^\beta), \ \ t\in [t_{n-1},t_n],
$$
and $\tilde{v}_N(t)=v_N(t-h)$. Taking the convergent subsequence in
the weak star topology in $H$, we know that the limit $v\in H$
satisfies
$$
C_\beta v^{\beta-1}\partial_tv^\beta =\Delta v
$$
in the distributional sense. The latter is the weak form of the
equation (\ref{pm2}).

Then we have proven the following result.

\begin{Thm}\label{Porous} Assume that $\Omega$ is a regular domain in the plane $R^2$. Given any $T>0$ and $m>0$.
Assume that the initial data $v_0=u_0^{1/m}\in H^1(\Omega)$. Let
$\alpha=1/m=2\beta-1$ with $\beta>1$. Then there is at least one
weak solution $v\in L^\infty H^1(\Omega)$ to the porous media
equation (\ref{pm2}).
\end{Thm}

Recall here that a mapping $v:[0,T]\to H^1(\Omega)$ is said \emph{a
weak solution} to (\ref{pm2}) if $v\in \in L^\infty H^1(\Omega)$
satisfies the evolution equation (\ref{pm2}) in the distributional
sense and $\lim_{t\to 0} u(t)=u_0$ in $L^2(\Omega)$ and
$u(t,x)=u_0(x)$ for $x\in \partial\Omega$ in the trace sense. The
latter will be simply said that $u(t)$ has the initial data $u_0$.

\section{$H^1$ weak solution to Ricci flow on $(S,g_0)$ with
symmetry}\label{sect3}

We now choose $\rho=\frac{2}{\fint_S e^{2u}dA}$ and fix any positive
constant $T>0$. Fix any $u_0\in H^1$ with $u(x)=u(-x)$, which is
symmetric about $r=\pi/2$. Then the Ricci flow equation under
consideration is
\begin{equation}\label{disc}
e^u\partial_t e^u=\Delta u-1+\frac{e^{2u}}{\fint e^{2u}dA}, \ \ in \
\ S\times (0,T]
\end{equation}
with the initial data $u(0)=u_0$. The symmetry condition is used to
get the compactness required for the minimization process below. We
now introduce a new concept of $H^1$ weak solution to the Ricci flow
(\ref{disc}).

\begin{Def}\label{weak}
We say $u(t)\in H^1$, $t\in [0,T]$ is a weak solution to the Ricci
flow (\ref{disc}) with the initial data $u_0$ if it satisfies
(\ref{disc}) in the sense of distribution and with the bounds
$$ \sup_{t\in [0,T]}\fint_S|\nabla u|^2+2udA-\log\fint_Se^{2u}dA\leq
\fint_S|\nabla u_0|^2+2u_0dA-\log\fint_Se^{2u_0}dA.
$$ and
$$
\int_{S\times [0,T]}|\partial_t e^{u(t)}|^2\leq
\frac{1}{2}\int_S|\nabla
u_0|^2+2u_0dA-\frac{1}{2}\log\fint_Se^{2u_0}dA.
$$
\end{Def}

In below, we may assume that $u_0$ is smooth (i.e., $u_0\in C^2$)
and normalize it such that the average $\bar{u}_0=\fint_S
u=\frac{1}{|S|}\int_Su_0dv_{g_0}$, otherwise, we choose an smoothly
approximation of $u_0$ and passing to subsequence limit.

For any $N>1$ be a large integer and for any $T>0$,let
$$
h=T/N, \ \ t_n=nh, \ \ n=0,1,2,...,N.
$$

Assume that we have constructed $u_j\in H^1(S)$, $0\leq j\leq n-1$
and $u_{n-1}$ is a minimizer of the functional
$$
J_{n-1}(u)=\frac{1}{2h}\fint_S
|e^u-e^{u_{n-2}}|^2dA+\frac{1}{2}\fint_S |\nabla
u|^2dA-\frac{1}{2}\log\fint_Se^{2u}dA
$$
over $\mathbf{A}=\{u\in H^1;\bar{u}=0\}$ and
 with the uniform bound
$$
\fint_S |\nabla u_{n-1}|^2dA-\log\fint_Se^{2u}dA\leq C.
$$
We want to get another $u_n$ with same bound by using the discrete
Morse method due to Rothe .
 Define the functional
$$
J_n(u)=\frac{1}{2h}\fint_S |e^u-e^{u_{n-1}}|^2dA+\frac{1}{2}\fint_S
|\nabla u|^2dA-\frac{1}{2}\log\fint_Se^{2u}dA
$$
in $ \mathbf{A}$. If $u_n$ is a minimizer of the functional on
$\mathbf{A}$, then it is clear that its Euler-Lagrange equation of
this functional is
\begin{equation}\label{EL}
e^{u_n}\frac{1}{h}(e^{u_n}-e^{u_{n-1}})=\Delta
u_n-\lambda_n+\frac{e^{2u_n}}{\fint_Se^{2u_n}}, \ \ on \ \ S,
\end{equation}
where
$$\lambda_n=1-\fint_Se^{u_n}\frac{1}{h}(e^{u_n}-e^{u_{n-1}}),$$ and
the minimizer satisfies the following estimate
\begin{equation}\label{est}
J_n(u_n)\leq J_n(u_{n-1})=\frac{1}{2}\int_S |\nabla
u_{n-1}|^2dA-\frac{1}{2}\log\fint_Se^{2u_{n-1}}dA.
\end{equation}

We want to minimize this functional on $\mathbf{A}$. First we need
to know that the infimum is finite. This follows from the fact that
after using the inequality $(\ref{chen})$, the leading term in the
functional $J_n$ is $\frac{1}{2h}\int_S |e^u-e^{u_{n-1}}|^2dA$.

 To get the minimizer, we only need to show that the minimizing
sequence is bounded in $\mathbf{A}$. Right from the the relation
$$
\inf_{H^1(S)} J_n(u)\leq J_n(u_{n-1})\leq C,
$$
 we know that the minimizing sequence $(u=u_{nk})$ satisfies
that
$$
\frac{1}{2h}\fint_S |e^u-e^{u_{n-1}}|^2dA+\frac{1}{2}\fint_S |\nabla
u|^2dA-\frac{1}{2}\log\fint_Se^{2u}dA\leq C.
$$
Using Chen-Li's Moser type inequality we know that the minimizing
sequence is uniformly bounded in $H^1$. By this we can pass to
subsequence to get a weakly convergent subsequence with its limit
$u_n$ as the minimizer.

Once this is done, we inductively define $u_N(t)\in H^1$ for $t\in
[-h,T]$ such that for $n=1,...,N$,
$$
u_N(t)=u_n, \  \ \lambda_N(t)=\lambda_n,
$$
on $[t_{n-1},t_n]$ and $u_N(t)=u_0$ on $[-h,0]$. We define, for
$1\leq n\leq N$,
$$
\partial_t e^{u_N(t)}=\frac{1}{h}(e^{u_n}-e^{u_{n-1}}), \ \ t\in [t_{n-1},t_n],
$$
and $\tilde{u}_N(t)=u_N(t-h)$. The estimate (\ref{est}) gives us
that $$ \sup_{t\in [0,T]}\fint_S|\nabla
u_n|^2dA-\frac{1}{2}\log\fint_Se^{2u_n}dA\leq \int_S|\nabla
u_0|^2dA-\frac{1}{2}\log\fint_Se^{2u_0}dA.
$$ and
$$
\int_{S\times [0,T]}|\partial_t e^{u_N(t)}|^2\leq
\frac{1}{2}\fint_S|\nabla
u_0|^2dA-\frac{1}{2}\log\fint_Se^{2u_0}dA+C_0.
$$

 The uniform upper bound of  $\int_S |\nabla
u_N|^2dA$ follows from the uniform bound of $\int_S e^{u_N}$ and the
assumption that $u(x)=u(-x)$ where the reflection is about
$r=\pi/2$. With the latter symmetry assumption we have the better
inequality in the sense of Moser that
$$
\fint_S e^{2u} dA\leq C_2exp\{\frac{1}{2}\fint_S|\nabla
u|^2dA+\fint_S 2u dA\}
$$
which is proved by W.Chen \cite{C} (see theorem II there). It is
here that the symmetry condition plays the role.

Note that $\{u_N\}$ is bounded in $L^\infty_tH^1$ and
$|\partial_te^{u_N}|$ is bounded in $L^2([0,t]\times S)$ for any
$t>0$. Then we may pass to weakly convergent subsequence of $u_N(t)$
in (\ref{EL}) to get the weak limit $u(t)\in L^\infty_tH^1(S)\cap
L^2([0,t]\times S)$ such that
\begin{equation}\label{mRF}
e^u\partial_t e^u=\Delta
u-\lambda(t)+\frac{e^{2u}}{\fint_Se^{2u}dA},
\end{equation}
where
$$
\lambda(t)=\lim \lambda_N(t)=1-\fint_Se^{u}\partial_te^u.
$$
The latter is called the relaxation parameter of the ricci flow and
the open problem is to prove that it decays to zero as time going to
infinity.

 Hence we have proven the following assertion.
\begin{Thm}\label{ingo}
For any positive constant $T>0$ and any $u_0\in H^1$ with
$u(x)=u(-x)$ with respect to $r=\pi/2$, there exists at least one
weak solution to (\ref{disc}) on $S\times [0,T]$ in the sense of the
weak form (\ref{mRF}).
\end{Thm}

Our argument above also work in the case of the standard sphere. It
may be used to define the weak solution to the Kaehler-Ricci flow.
It is an open question if we have stability result or uniqueness
result for the weak solution above. We remark that it is possible to
get a little more regularity of the weak solution for smooth initial
data $u_0$. For higher dimensions, except the concept of Ricci flow
with surgery, it is not clear if one can define the weak Ricci flow.

\section{Another discrete Morse flow}\label{sect4}
Fix $\lambda\in (0,1)$. This parameter will play the regularization
role in the discrete Ricci flow below. We study the following
modified Ricci flow
\begin{equation}\label{disc2}
e^u\partial_t e^u=\Delta u-\lambda+\frac{e^{2u}}{\fint e^{2u}dA}, \
\ in \ \ S\times (0,T]
\end{equation}
with the initial data $u(0)=u_0\in H^1$. In this case, the variation
structure related to the right side of (\ref{disc2}) is
$$
I(u)=\frac{1}{2}\fint_S(|\nabla u|^2+2\lambda
u)dA-\frac{1}{2}\log\fint_S e^{2u}dA.
$$
The advantage of this function is that the inequality (\ref{chen1})
implies that
\begin{equation}\label{f1}
I(u)\geq 2)(\lambda-1)\fint_S udA
\end{equation}
and
\begin{equation}\label{f2}
I(u)\geq \frac{1-\lambda}{2}\fint |\nabla
u|^2dA-\frac{1-\lambda}{2}\log\fint_S e^{2u}dA.
\end{equation}
These two facts will help us to make the infimum of the functional
$\hat{J}$ (see below) in $H^1$ be finite.

As in last section, for any $N>1$ be a large integer and for any
$T>0$,let
$$
h=T/N, \ \ t_n=nh, \ \ n=0,1,2,...,N.
$$

Assume that we have constructed $\hat{u}_j\in H^1(S)$, $0\leq j\leq
n-1$ and $\hat{u}_{n-1}$ is a minimizer of the functional
$$
\hat{J}_{n-1}(u)=\frac{1}{2h}\fint_S
|e^u-e^{\hat{u}_{n-2}}|^2dA+\frac{1}{2}\fint_S (|\nabla
u|^2+2\lambda u)dA-\frac{1}{2}\log\fint_Se^{2u}dA
$$
over $H^1$ and
 with the uniform bound
$$
\fint_S (|\nabla u_{n-1}|^2+2\lambda u)dA-\log\fint_Se^{2u}dA\leq C.
$$

Let $$ \hat{J}(u)=\frac{1}{2h}\fint_S
|e^u-e^{\hat{u}_{n-1}}|^2dA+I(u)
$$
on $H^1$.
 Because of the inequality (\ref{chen1}) we can easily show as
before that the minimizing sequence of $\hat{J}$ is bounded in
$H^1$. In fact, if the infimum is $\infty$, then we may have
sequence $w_j\in H^1$ such that $\hat{J}(w_j)\to-\infty$. By
(\ref{f2}) or (\ref{f1}) we know that $\int_Se^{2w_j}\to\infty$ or
$\fint_Sw_j\to\infty$, which is impossible since in this case the
leading term is $$\frac{1}{h}\fint_S|e^{w_j}-e^{u_{n-1}}|^2.$$

Hence the infimum is finite. By this we can pass to subsequence to
get a weakly convergent subsequence with its limit $\hat{u}_n$ as
the minimizer of $\hat{J}$. The Euler-lagrange equation of $\hat{J}$
is
\begin{equation}\label{EL2}
e^{\hat{u}_n}\frac{1}{h}(\hat{e}^{u_n}-e^{\hat{u}_{n-1}})=\Delta
\hat{u}_n-\lambda+\frac{e^{2\hat{u}_n}}{\fint_Se^{2\hat{u}_n}}, \ \
on \ \ S.
\end{equation}
By the standard elliptic regularity theory we know that $\hat{u}_n$
is smooth.

Then we inductively define $\hat{u}_N(t)\in H^1$ for $t\in [-h,T]$
such that for $n=1,...,N$,
$$
\hat{u}_N(t)=\hat{u}_n,
$$
on $[t_{n-1},t_n]$ and $\hat{u}_N(t)=u_0$ on $[-h,0]$. We define,
for $1\leq n\leq N$,
$$
\partial_t e^{\hat{u}_N(t)}=\frac{1}{h}(e^{\hat{u}_n}-e^{\hat{u}_{n-1}}), \ \ t\in [t_{n-1},t_n]
$$
 the minimizing property of $\hat{u}_n$ gives
us the desired energy bound to obtain a convergent subsequence
$\hat{u}_N$ and the limit $\hat{u}\in L^\infty_tH^1(S)\cap
L^2([0,t]\times S)$ such that
\begin{equation}\label{mRF2}
e^u\partial_t e^u=\Delta u-\lambda+\frac{e^{2u}}{\fint_Se^{2u}dA}.
\end{equation}

Hence we have the following result.

\begin{Thm}\label{ingo2}
For any positive constant $T>0$ and any $u_0\in H^1$, there exists
at least one weak solution to (\ref{disc2}) on $S\times [0,T]$ in
the sense of the weak form (\ref{mRF2}).
\end{Thm}

Here the weak solution of (\ref{disc2}) is defined similar to the
definition \ref{weak}.

{\bf Acknowledgement}. The paper has been done when the fist named
author is visiting Math.Inst., Gottingen University, Germany, in
February 2012 and he would like to thank the host institute for
hospitality.

\end{document}